\documentclass[twocolumn,showpacs,preprintnumbers,amsmath,amssymb]{revtex4}

\usepackage{dcolumn}
\usepackage{bm}
\usepackage{amssymb}
\usepackage{amsthm,epsf}
\usepackage{amsfonts}
\usepackage{epsfig}
\usepackage{graphicx}
\usepackage{latexsym}

\begin{document}
\title{\Large\bf Degree-distribution stability of scale-free networks}

\author{Zhenting Hou$^{1\ast}$}
\author{Xiangxing Kong$^{1\star}$}
\author{Dinghua Shi$^{1,2\dagger}$}
\author{Guanrong Chen$^{3\ddagger}$}

\affiliation{$^{1}$School of Mathematics,
      Central South University, Changsha 410083, China\\
$^{2}$Department of Mathematics, Shanghai University, Shanghai 200444, China\\
$^{3}$Department of Electronic Engineering, City University of
      Hong Kong, Hong Kong, China} 

\medskip

\date{\today}

\begin{abstract}
Based on the concept and techniques of first-passage probability
in Markov chain theory, this letter provides a rigorous proof for
the existence of the steady-state degree distribution of the
scale-free network generated by the Barab\'asi-Albert (BA) model,
and mathematically re-derives the exact analytic formulas of the
distribution. The approach developed here is quite general,
applicable to many other scale-free types of complex networks.
\end{abstract}



\maketitle

{\it Introduction}. The intensive study of complex networks is
pervading all kinds of sciences today, ranging from physical to
biological, even to social sciences. Its impact on modern
engineering and technology is prominent and will be far-reaching.
Typical complex networks include the Internet, the World Wide Web,
wired and wireless communication networks, power grids, biological
neural networks, social relationship networks, scientific
cooperation and citation networks, and so on. Research on
fundamental properties and dynamical features of such complex
networks has become overwhelming.

In the investigation of various complex networks, the degree
distributions are always the main concerns because they
characterize the fundamental topological properties of the
underlying networks.

Noticeably, for a ring-shape regular graph$^{[1]}$ of whatever
size, where every vertex is connected to its $K$
nearest-neighboring vertices, all vertices have the same degree
$K$. For the well-known Erd\"os-R\'enyi random graph model$^{[2]}$
with $n$ vertices and $m$ edges, the degree distribution of
vertices is approximately Poisson with mean value $2m/n$. For the
small-world network proposed by Watts and Strogatz$^{[1]}$, the
degree distribution of vertices also follows Poisson distribution
approximately.

A common feature of the above models is that the degree
distribution of vertices has a characteristic size $\langle
k\rangle$. In contrast, Barab\'asi and Albert$^{[3]}$ found that
for many real-world complex networks, e.g., the WWW, the fraction
$P(k)$ of vertices with degree $k$ is proportional over a large
range to a ``scale-free'' power-law tail: $k^{-\gamma}$, where
$\gamma$ is a constant independent of the size of the network.
Thus, the fraction $P(k)$ of vertices with degree $k$ is referred
to as the degree distribution of a scale-free network. To explain
this phenomenon, they proposed the following network-generating
mechanism$^{[3]}$, known as the BA model:

``$\cdots$ starting with a small number ($m_0$) of vertices, at
every time step we add a new vertex with $m$ ($\le m_0$) edges
that link the new vertex to $m$ different vertices already present
in the system. To incorporate preferential attachment, we assume
that the probability $\Pi$ that a new vertex will be connected to
a vertex depends on the connectivity $k_i$ of that vertex, so that
$\Pi(k_i)=k_i/\sum_j k_j$. After $t$ steps the model leads to a
random network with $t+m_0$ vertices and $mt$ edges.''
\medskip

In [3], computer simulation showed that for the BA model the
degree distribution of the network has a power law form with the
exponent $\gamma=2.9\pm0.1$. In [4], a heuristic argument based on
the mean-field theory led to an analytic solution
$P(k)\,\sim\,2m^2k^{-3}$, namely $\gamma=3$. To derive the
following dynamic equation:
\[
\frac{\partial k_i}{\partial t}=m\Pi(k_i)=\frac{k_i}{2t},\qquad
k_i(i)=m\,,
\]
it was assumed$^{[4]}$ that the probability (can be interpreted as
a continuous rate of change of $k_i$) for an existing vertex with
degree $k_i$ to receive a new connection from the new vertex is
exactly equal to $m\Pi(k_i)$, which is simultaneously proportional
to both the degree $k_i(t)$ of the existing vertex $i$ and the
number $m$ of the new edges that the new vertex brings in, at time
$t$. For notational convenience, this assumption will be simply
referred to as the ``$m\Pi$-hypothesis'' in this paper.

In all the consequent works related to the BA model, this
$m\Pi$-hypothesis plays a fundamental role. For example, Krapivsky
et al.$^{[5]}$ replaced the degree $k_i(t)$ of vertex $i$ at time
$t$ by the total number $N_k(t)$ of degree-$k$ vertices over the
whole network at time $t$, thereby obtaining its rate equation
\[
\frac{d N_k(t)}{dt}=m\,\frac{(k-1)N_{k-1}(t)-kN_k(t)}{\sum_k k
N_k(t)}+\delta_{km}\,,
\]
where $\delta_{km}$ accounts for new vertices bringing in new
edges. In this study, the $m\Pi$-hypothesis was adopted in the
derivations. Assuming that the steady-state degree distribution
exists, using the law of large numbers ($\frac{N_k(t)}{t}\to P(k)$
as $t\to\infty$), they showed that the difference equation of
$P(k)$ has an analytic solution
\begin{eqnarray*}
P(k)=\frac{4}{k(k+1)(k+2)}
\end{eqnarray*}
for the BA model with $m=1$. They also pointed out that only the
linear preferential attachment scheme can lead to the scale-free
structure but any nonlinear one will not.

Dorogovtsev et al.$^{[6]}$ considered $k_i(t)$ as a random
variable and defined $P(k,i,t)$ to be the probability that vertex
$i$ has exactly $k$ edges at time $t$, where vertex $i$ is the
vertex that was being added to the network at time $t=i$,
$i=1,2,\cdots\,$. Moreover, they used the average of all vertex
degrees as the network degree: $P(k,t)\triangleq
\frac{1}{t}\sum_{i=1}^t P(k,i,t)$. They introduced a more general
attraction model and allowed multiple edges between vertices,
where each new vertex has an initial attraction degree $A$.
Simultaneously, $m$ new directed edges coming out from
non-specified vertices are introduced with the probability $\Pi$,
therefore $k=A+q$ with $q$ being the in-degree of vertices.
Consequently, when every new vertex is the source of the $m$ new
edges like in the BA model, the attraction model makes more sense
than the BA model. They first arrived at the master equation of
$P(k,i,t)$ and then by summing all $i$'s together they were able
to derive the following equation:
\begin{eqnarray*}
P(k,t+1)&=&\frac{k-1}{2t}\,P(k-1,t)+
\left(1-\frac{k}{2t}\right)P(k,t)\\
& &+\delta_{mk}+O\left(\frac{P(k,t)}{t}\right).
\end{eqnarray*}
To that end, by assuming the existence of $P(k)$ [note that
actually an additional assumption of $\lim_{t\to\infty}\,
t[P(k,t+1)-P(k,t)]=0$ is also needed], they obtained a difference
equation for $P(k)$. Finally, solving the equation gave an
analytic solution
\begin{eqnarray*}
P(k)=\frac{2m(m+1)}{k(k+1)(k+2)}.
\end{eqnarray*}
Here, it should be pointed out that if multiple edges are not
allowed, then the $m\Pi$-hypothesis is still needed.

As a side note, Dorogovtsev et al.$^{[7]}$ also considered the
effect of accelerating growth, which is proportional to the power
of the time variable $t$ at each time step. However, this destroys
the scale-free feature and degree-distribution stability of the
network.

\medskip
Afterwards, Bollob\'as$^{[8]}$ made a general comment on the BA
model:

``From a mathematical point of view, however, the description
above, repeated in many papers, does not make sense. The first
problem is getting started. The second problem is with the
preferential attachment rule itself, and arises only for $m\geq2$.
In order to prove results about the BA model, one must first
decide on the details of the model itself. It turns out to be
convenient to allow multiple edges and loops.''

Consequently, he and his coauthors recommended a so-called LCD
model, as follows:

``We start with the case $m=1$. Consider a fixed sequence of
vertices $v_1,\,v_2\,,\cdots\,$. We shall inductively define a
random graph process $\{G_1^t\}_{t\geq0}$ so that $G_1^t$ is a
directed graph on $\{v_i:\,1\le i\le t\}$, as follows. Start with
$G_1^0$ the ``graph'' with no vertices, or with $G_1^1$ the graph
with one vertex and one loop. Given $G_1^{t-1}$, form $G_1^t$ by
adding the vertex $v_t$ together with a single edge directed from
$v_t$ to $v_i$, where $i$ is chosen randomly with
\[
\Pi(i=s)=\left\{ \begin{array}{ll}
 \frac{d_{G_1^{t-1}(v_s)}}{2t-1}, & 1\le s\le t-1\\[5pt]
 \frac{1}{2t-1}, & s=t.
\end{array}
\right.
\]
For $m>1$ we define the process $\{G_m^t\}_{t\geq0}$ by running
the process $\{G_1^t\}$ on a sequence
$v_1^{\prime},\,v_2^{\prime},\,\cdots\,$; the graph $G_m^t$ is
formed from $G_1^{mt}$ by identifying the vertices
$v_1^{\prime},\,v_2^{\prime},\,\cdots\,,v_m^{\prime}$ to form
$v_1$, identifying
$v_{m+1}^{\prime},\,v_{m+2}^{\prime},\,\cdots\,,v_{2m}^{\prime}$
to form $v_2$, and so on.''

For graph $G_m^n$, let $\#_m^n(d)$ be the number of vertices of
$G_m^n$ with in-degree equal to $d$, i.e., with (total) degree
$m+d$, and set
\[a_{m,d}=\frac{2m(m+1)}{(d+m)(d+m+1)(d+m+2)}.\]
Bollob\'as et al.$^{[9]}$ rigorously proved the following result:
\[\lim_{n\rightarrow\infty}E[\#_m^n(d)]/n=a_{m,d}.\]
Then, based on the martingale theory, they proved that
$\#_m^n(d)/n$ converges to $a_{m,d}$ in probability.
\medskip

It has been observed that most real-world and simulated networks
follow certain rules to add or remove their vertices and edges,
which are not entirely random. More importantly, at each time
step, these rules are applied only to the previously formed
network, therefore the process has prominent Markovian properties.
Shi et al.$^{[10]}$ established a close relationship between the
BA model and Markov chains. According to the evolution of the BA
model, the degree $k_i(t)$ of vertex $i$ at time $t$ constitutes a
nonhomogeneous Markov chain as time evolves. Thus, all vertices
together form a family of Markov chains. Consequently, based on
the Markov chain theory, starting from an initial distribution and
iteratively multiplying the state-transition probability matrices,
the final network degree distribution can be easily obtained.
Lately, Shi et al.$^{[11]}$ developed an evolving network model by
using an anti-preferential attachment mechanism, which can
generate scale-free networks with power-law exponents varying
between $1 \thicksim 4$. There are several modified and
generalized BA models in the literature, including such as the
local-world BA model$^{[12]}$, which will not be listed and
reviewed here.

All in all, the BA model indeed is a breakthrough discovery with
significant impact on network science today. Therefore, it is
quite important to support the model with a rigorous mathematical
foundation.

It is clear from the above discussions that two key questions need
to be carefully answered for the BA model: 1) For the case of
$m\geq2$, can one find a scheme of adding new edges from the new
vertex to the existing ones that has a probability precisely equal
to $m\Pi$? This is the key of the BA modeling. 2) Does the
steady-state degree distribution of the network exist and, if so,
what is it? This is the key to the validity of the mean-field,
rate-equation, master-equation, and Markov-chain approaches. The
present paper will give complete answers to these two questions.

\bigskip

{\it Degree-distribution stability}. To start, consider the first
question. Recall that Holme and Kim$^{[13]}$ proposed a scheme for
new edge connection: When a new vertex comes into the network, the
first edge connects to an existing vertex with the preferential
attachment probability $\Pi$. After that, the rest $m-1$ edges
randomly connect with probability $p$ to the vertices in the
neighborhood of the vertex that the first edge was connected to,
or connect with probability $1-p$ to those vertices that the first
edge did not connect to. Here, consider this approach with $p=1$
in the following scenario: When a new vertex comes into the
network, the first edge connects to an existing vertex with the
specified preferential attachment probability $\Pi$, same as
above. Yet, the rest $m-1$ edges simultaneously connect to $m-1$
vertices randomly chosen from inside the neighborhood of the
vertex that the first edge was connected to. By random sampling
theory this is equivalent to the above Holme-Kim scheme which
continually connects the edges to $m-1$ vertices randomly chosen
from inside the neighborhood without allowing multiple edges. For
this special scheme, the following result can be rigorously
proved.
\medskip

\noindent{\bf Proposition.} For the BA model with the above
special attachment scheme, if vertex $i$ has degree $k_i(t)$ at
time $t$, then the probability that vertex $i$ receives a new edge
from the new vertex at time $t+1$ is exactly equal to $m\Pi(k_i)$.
\medskip

\noindent{\it Proof.} Let $P_i(t+1)$ be the probability of vertex
$i$ receiving a new edge from vertex $t+1$ at time $t+1$. Then,
\begin{eqnarray*}
P_i(t+1)&=&\frac{k_i(t)}{\sum_j k_j(t)}+\sum_{l\in O_i(t)}
\frac{k_l(t)}{\sum_j k_j(t)}\,
\frac{C_{k_l(t)-1}^{m-2}}{C_{k_l(t)}^{m-1}}\\
&=&\frac{k_i(t)}{\sum_j k_j(t)}+\sum_{l\in
O_i(t)}\frac{m-1}{\sum_j k_j(t)} = m\,\frac{k_i(t)}{\sum_j
k_j(t)}\,,
\end{eqnarray*}
where
\begin{small}\[ \frac{C_{k_l(t)-1}^{m-2}}{C_{k_l(t)}^{m-1}}
=\frac{(k_l(t)-1)!/[(m-2)!(k_l(t)-m+1)!]}{k_l(t)!/[(m-1)!(k_l(t)-m+1)!]}
=\frac{m-1}{k_l(t)},
\]\end{small}

\noindent which is the probability of choosing vertex $i$, among
the $m-1$ vertices that were randomly chosen from inside the
neighborhood $O_l(t)$ of vertex $l$, to perform simultaneous
connections.

The Proposition answers the first question posted above and shows
that the special Holme-Kim preferential attachment scheme is one
way to implement the $m\Pi$-hypothesis.
\medskip

In order to prove the degree-distribution stability of the general
BA network, the BA model is specified first. Start with a complete
graph with $m_0$ vertices, which has a total degree
$N_0=m_0(m_0-1)$, and denote these vertices by $-m_0,\cdots,-1$,
respectively. In all the following derivations, the
$m\Pi$-hypothesis will be assumed. The general BA networks will be
further discussed in the last section below.

Following Dorogovtsev et al.$^{[6]}$, consider the degree $k_i(t)$
as a random variable, and let $P(k,i,t)=P\{k_i(t)=k\}$ be the
probability of vertex $i$ having degree $k$ at time $t$, and
moreover let the network degree distribution be the average over
all its vertices at time $t$, namely,
\begin{eqnarray*}
P(k,t)\triangleq \frac{1}{t+m_0}\sum_{i=-m_{0},i\neq0}^t P(k,i,t).
\end{eqnarray*}

Recall that $k_i(t)$ is a random variable for any fixed $t$ and it
is a nonhomogeneous Markov chain for variable $t$$^{[10]}$. Under
the $m\Pi$-hypothesis, the state-transition probability of this
Markov chain is given by
\begin{equation}
P\{k_i(t+1)=l\,|\,k_i(t)=k\}= \left\{
\begin{array}{ll}
1-\frac{k}{2t+\frac{N_0}{m}},~l=k\\
\frac{k}{2t+\frac{N_0}{m}},~l=k+1\\
0,~~~~{\rm otherwise},
\end{array}
\right.
\end{equation}
where $k=1,2,\cdots,m+t-i$, and $i=1,2,\cdots\,$.

The existence of the steady-state degree distribution for this
specified BA network can be proved in three steps as follows.
Detailed derivations are supplied in the Appendix of the paper.
\medskip

\noindent{\bf 1.} Consider the first-passage probability of the
Markov chain:
\[
f(k,i,t)=P\{k_i(t)=k,\,k_i(l)\not=k,\,l=1,2,\cdots,t-1\}.
\]
Then, the relationship between the first-passage probability and
the vertex degrees is established.\medskip

\noindent{\bf Lemma 1.} Under the $m\Pi$-hypothesis, for the BA
model with $k>m$,
\begin{eqnarray}
f(k,i,s)&\!\!=\!\!&P(k-1,i,s-1)\,\frac{k-1}{2(s-1)+\frac{N_0}{m}}\,,\\[5pt]
P(k,i,t)&\!\!=\!\!&\sum_{s=i+k-m}^tf(k,i,s)\,\prod_{j=s}^{t-1}
           \left(1-\frac{k}{2j+\frac{N_0}{m}}\right).
\end{eqnarray}

\noindent{\bf 2.} Under the $m\Pi$-hypothesis, using the
state-transition probability of the Markov chain, one first finds
the expression of $P(m,t)$, as follows:
\begin{small}\begin{eqnarray*}
P(m,t)&=&\prod\limits_{i=1}^{t-1}\left(1-\frac{m}{2i+\frac{N_{0}}{m}}\right)
\frac{i+m_0}{i+1+m_0}\\
&&\times\left[\,P(m,1)+\sum\limits_{l=1}^{t-1}
\frac{\frac{1}{l+1+m_0}}{\prod\limits_{j=1}^{l}
\left(1-\frac{m}{2j+\frac{N_{0}}{m}}\right)\frac{j+m_0}{j+1+m_0}}\right]\\
&=&\frac{1}{t+m_0}\prod\limits_{i=1}^{t-1}
\left(1-\frac{m}{2i+\frac{N_{0}}{m}}\right)\\
&&\times\left[(1+m_0)P(m,1)
+\sum\limits_{l=1}^{t-1}\prod\limits_{j=1}^{l}
\left(1-\frac{m}{2j+\frac{N_{0}}{m}}\right)^{-1}\right].
\end{eqnarray*}\end{small}

Then, one can show the existence of the limit
$\lim_{t\to\infty}P(m,t)$ by using the following classical
Stolz-Ces\'aro Theorem in Calculus.\medskip

\noindent{\bf Stolz-Ces\'aro Theorem$^{[14]}$.} In sequence
$\{\frac{x_n}{y_n}\}$, assume that $\{y_n\}$ is a monotone
increasing sequence with $y_n\rightarrow\infty$. If the limit
$\lim\limits_{n\rightarrow\infty}\frac{x_{n+1}-x_n}{y_{n+1}-y_n}=l$
exists, where $-\infty\leq l\leq+\infty$, then
$\lim\limits_{n\rightarrow\infty}\frac{x_n}{y_n}=l$.\medskip

\noindent{\bf Lemma 2.}\ Under the $m\Pi$-hypothesis, for the BA
model, the limit\ $\lim_{t\to\infty}P(m,t)$ exists and is
independent of the initial network:
\begin{equation}
P(m)\triangleq \lim_{t\to\infty}P(m,t)=\frac{2}{m+2}>0\,.
\end{equation}

\noindent{\bf 3.} Under the $m\Pi$-hypothesis, similarly, one
finds the expression of $P(k,t)$ using the first-passage
probability of the Markov chain, and then shows the existence of
the limit $\lim_{t\to\infty}P(k,t)$ by using the Stolz-Ces\'aro
Theorem, if the limit\ $\lim_{t\to\infty}P(k-1,t)$ exists.
\medskip

\noindent{\bf Lemma 3.} Under the $m\Pi$-hypothesis, for the BA
model with $k>m$, if the limit\ $\lim_{t\to\infty}P(k-1,t)$ exists
then the limit\ $\lim_{t\to\infty}P(k,t)$ also exists:
\begin{equation}
P(k)\triangleq
\lim_{t\to\infty}P(k,t)=\frac{k-1}{k+2}\,P(k-1)>0\,.
\end{equation}

Finally, by mathematical induction, it follows from Lemmas 2 and 3
that the steady-state degree distribution of the specified BA
network exists. To this end, by solving the difference equation
(5) iteratively, one arrives at the following conclusion.
\medskip

\noindent{\bf Theorem 1.} Under the $m\Pi$-hypothesis, for the BA
model with $k\geq m$, the steady-state degree distribution exists,
independent of the initial network, and is given by
\begin{equation}
P(k)=\frac{2m(m+1)}{k(k+1)(k+2)}\ \thicksim\ 2m^2k^{-3}>0\,.
\end{equation}

Clearly, this degree distribution formula is consistent with the
formula obtained by Dorogovtsev et al.$^{[6]}$ and Bollob\'as et
al.$^{[9]}$, which allow multiple edges and loops.
\bigskip

{\it Discussion}. Bollob\'as$^{[8]}$ once discussed the BA
description (the $m\Pi$-hypothesis) of preferential attachment in
detail. His result gives a range of models fitting the BA
description with very different properties. When $m\geq2$, as a
new vertex comes in, it is no problem for its first edge to
preferentially connect to an existing vertex. But what about the
other $m-1$ new edges? This question was not carefully addressed
before. Clearly, after the first edge has been connected from the
new vertex to an existing vertex, the preferential attachment
probability $\Pi$ is no longer the same if later operations do not
allow multiple edges and loops. It is also clear that when
$m\geq2$, the probability of vertex $i$ receiving a new edge is
always greater than $\Pi$. But what is it? On the other hand, it
is also possible that the probability of vertex $i$ receiving a
new edge depends on other vertex degrees. Barab\'asi always
emphasizes the $m\Pi$-hypothesis but did not discuss this ``how''
question either. Thus, two questions arise: 1) For the BA model,
or for any other BA-like model, how to prove the
degree-distribution stability if the $m\Pi$-hypothesis holds only
approximately? 2) Is there a preferential attachment scheme for
$m\geq2$ such that the probability of vertex $i$ receiving a new
edge is independent of other vertex degrees?

To answer these two questions, a new preferential attachment
scheme is proposed and discussed in [15], where a new vertex will
be simultaneously connected to $m$ different vertices and it is
assumed that the preferential attachment probability $\Pi$ is
proportional to the sum of the degrees $k_{i_1},\cdots,k_{i_m}$ of
those vertices. They showed that the probability that the existing
vertex $i$ received an edge from the new vertex is independent of
other vertex degrees, namely,
\begin{eqnarray*}
\Pi_m^{t+1}(k_{i}(t))=\frac{m_0+t-m}{m_0+t-1}
\frac{k_i(t)}{2mt+N_0}+\frac{m-1}{m_0+t-1},
\end{eqnarray*}
where $m_0$ is the number of vertices and $N_0$ is the total
degree in the initial network. Consequently, under the
$(a_tk_i(t)+b_t)(1+o(1)_{k_i(t),t})$-hypothesis and some mild
conditions, they proved the degree-distribution stability of
Barab\'asi-Albert type networks. Especially, the power-law
exponent of the network degree distribution in this new
preferential attachment scheme is $\gamma=2m+1$.

Finally, it should be emphasized that the theory and scheme
developed in this paper has great generality$^{[16]}$, in the
sense that it can be applied to many BA-like modified and
generalized models, such as the LCD model of Bollob\'as et
al.$^{[9]}$, the attraction model of Dorogovtsev et al.$^{[6]}$,
the local-world BA-like model of Li and Chen$^{[12]}$, and the
evolving network model of Shi et al.$^{[11]}$, etc.

We summarize the results and findings in this paper as follows:
(1) Our proving method differs from the one based on martingale
theory, and can be applied to many other scale-free types of
complex networks; (2) We do not need to change the BA model, e.g.,
to allow multiple edges and loops; (3) We provide a special
Holme-Kim preferential attachment scheme such that the
``$m\Pi$-hypothesis'' holds.

This research was supported by the National Natural Science
Foundation under Grant No. 10671212, and by the NSFC-HKRGC Joint
Research Projects under Grant N-CityU107/07.

$\ast$ Email address: zthou@csu.edu.cn

$\star$ Email address: kongxiangxing2008@163.com

$\dagger$ Email address: shidh2001@263.net

$\ddagger$ Email address: gchen@ee.cityu.edu.hk


\begin{thebibliography}{99}

\bibitem{W}
Watts D. J. and Strogatz S. H., {\it Nature} {\bf 393}, 1998,
440-442

\bibitem{E}
Erd\"os P. and R\'enyi A., {\it Publications Mathematicae} {\bf
6}, 1959, 290-297

\bibitem{B1}
Barab\'asi A.-L. and Albert R., {\it Science} {\bf 286}, 1999,
509-512

\bibitem{B2}
Barab$\acute{a}$si A.-L., Albert R. and Jeong H., {\it Physica A}
{\bf 272}, 173 (1999).

\bibitem{K}
Krapivsky P. L., Redner S. and Leyvraz F., {\it Phys. Rev. Lett.}
{\bf 85}, 2000, 4629-4632

\bibitem{D1}
Dorogovtsev S. N., Mendes J. F. F. and Samukhin A. N., {\it Phys.
Rev. Lett.} {\bf 85}, 2000, 4633-4636

\bibitem{D2}
Dorogovtsev S. N. and Mendes J. F. F., {\it Phys. Rev. E} {\bf
63}, 2001, 025101

\bibitem{Bo}
Bollob\'as B., {\it Handbook of Graphs and Networks: From the
Genome to the Internet} (Bornholdt S. and Schuster H. G. eds.),
Wiley-VCH, 2002, 1-34

\bibitem{Bo1}
Bollob\'as B., Riordan O. M., Spencer J. and Tusn\'ady G., {\it
Random Structures and Algorithms} {\bf 18}, 2001, 279-290

\bibitem{S1}
Shi D. H., Chen Q. H. and Liu L. M., {\it Phys. Rev. E} {\bf 71},
2005, 036140

\bibitem{S2}
Shi D. H., Liu L. M., Zhu X. and Zhou H. J., {\it Europhys. Lett.}
{\bf 76}, 2006, 731-737

\bibitem{L}
Li X. and Chen G. R., {\it Physica A} {\bf 328}, 2003, 274-286

\bibitem{H}
Holme P. and Kim B. J., {\it Phys. Rev. E} {\bf  65}, 2002, 026107

\bibitem{St}
Stolz O., {\it Vorlesungen uber allgemiene Arithmetic}, Teubner,
Leipzig 1886

\bibitem{Ho}
Hou Z. T. et al., On the degree-distribution stability of
Barab\'asi-Albert type networks, 2008, preprint.

\end{thebibliography}
\end{document}


\pagestyle{plain}
\title{{\bf Appendix to\\[8pt]
``Degree-distribution stability of scale-free networks''}}
\author{Zhenting Hou$^1$\  Xiangxing Kong$^1$\  Dinghua Shi$^{1,2}$\  Guanrong Chen$^3$\\[3pt]
{\small $^{1}$School of Mathematics, Central South University, Changsha 410083, China}\\[3pt]
{\small $^{2}$Department of Mathematics, Shanghai University, Shanghai 200444, China}\\[3pt]
{\small $^{3}$Department of Electronic Engineering, City University of Hong Kong, Hong Kong, China}}%

\date{\today}

\maketitle

To provide a rigorous proof of the degree-distribution stability
of the scale-free network generated by the BA model, some
parameters are specified as follows: (i) start with a complete
graph with $m_0$ vertices, which has the total degree
$N_0=m_0(m_0-1)$, and denote these vertexes by $-m_0,\cdots,-1$,
respectively; (ii) assume that at each time step $t$, the
probability of the new vertex connecting to an existing vertex $i$
is exactly equal to $m\Pi(k_i(t))$.

Here, in (ii), the preferential attachment probability is
simultaneously proportional to both the degree $k_i(t)$ of the
existing vertex $i$ and the number $m$ of new edges that the new
vertex brings in, at time $t$. For notational convenience, this
assumption will be referred to as the ``$m\Pi$-hypothesis'' below.

Observe that the degree $k_i(t)$ of vertex $i$ at time $t$ is a
random variable$^{[6]}$. Let $P(k,i,t)=P\{k_i(t)=k\}$ denote the
probability of vertex $i$ having degree $k$ at time $t$, and
define the degree distribution of the whole network by the average
value of probabilities of vertex degrees
\begin{equation}
P(k,t)\triangleq
\frac{1}{t+m_0}\,\sum_{i=-m_{0},i\neq0}^tP(k,i,t).
\end{equation}

Observe also that the degree $k_i(t)$ as a process in time $t$ is
an nonhomogeneous Markov chain$^{[10]}$. Thus, for
$k=1,2,\cdots,t+i-m$, the state transition probabilities of the
Markov chain, under the $m\Pi$-hypothesis, are given by
\begin{eqnarray}
P\{k_i(t+1)=l|k_i(t)=k\}=\left\{\begin{array}{ll}
1-\frac{k}{2t+\frac{N_0}{m}}, &l=k\\[5pt]
\frac{k}{2t+\frac{N_0}{m}}, &l=k+1\\[5pt]
0, &{\rm otherwise}.
\end{array}\right.
\end{eqnarray}

\section{The BA model with $m=1$}

\ \ \ \ Denote the first-passage probability of the Markov chain
by $f(k,i,t)=P\{k_i(t)=k,\ k_i(l)\neq k,\ l=1,2,\cdots,t-1\}$.
First, the relationship between the first-passage probability and
the vertex degrees is established.
\medskip

\textbf{Lemma 1} For $k>1$,
\begin{eqnarray}
f(k,i,s)&=&P(k-1,i,s-1)\,\frac{k-1}{2(s-1)+N_{0}},\\
P(k,i,t)&=&\sum\limits_{s=i+k-1}^t
           f(k,i,s)\,\prod\limits_{j=s}^{t-1}
           \left(1-\frac{k}{2j+N_{0}}\right).
\end{eqnarray}

\textbf{Proof} First, consider Eq. (3). The degree of a vertex is
always nondecreasing, and increasing at most by $1$ each time,
according to the construction rule of the BA model. Thus, it
follows from the Markovian properties that
\begin{eqnarray*}
& &f(k,i,s)=P\{k_i(s)=k,\ k_i(l)\neq k,\ l=1,2,\cdots,s-1\}\\
& &=P\{k_i(s)=k,\ k_i(s-1)=k-1,\ k_i(l)\neq k,\ l=1,2,\cdots,s-2\}\\
& &=P\{k_i(s)=k,\ k_i(s-1)=k-1\}\\
& &=P\{k_i(s-1)=k-1\}\,P\{k_i(s)=k|k_i(s-1)=k-1\}\\
& &=P(k-1,i,s-1)\,\frac{k-1}{2(s-1)+N_{0}}.
\end{eqnarray*}

Second, observe that the earliest time for the degree of vertex $i$
to reach $k$ is at step $k+i-1$, and the latest time to do so is at
step $t$. After this vertex degree becomes $k$, it will not increase
any more. Thus, Eq. (4) is proved.
\medskip

\textbf{Lemma 2} (Stolz-Ces\'aro Theorem) In sequence
$\{\frac{x_n}{y_n}\}$, assume that $\{y_n\}$ is a monotone
increasing sequence with $y_n\rightarrow\infty$. If the limit
$\lim\limits_{n\rightarrow\infty}\frac{x_{n+1}-x_n}{y_{n+1}-y_n}=l$
exists, where $-\infty\leq l\leq+\infty$, then
$\lim\limits_{n\rightarrow\infty}\frac{x_n}{y_n}=l$.

\textbf{Proof} This is a classical result, see [14].
\medskip

\textbf{Lemma 3} For the probability $P(k,t)$ defined in Eq. (1),
the limit $\lim\limits_{t\rightarrow\infty}P(1,t)$ exists and is
independent of the initial network; moreover,
\begin{eqnarray}
P(1)\triangleq\lim\limits_{t\rightarrow\infty}P(1,t)=\frac{2}{3}>0.
\end{eqnarray}

\textbf{Proof} From the construction of the BA model or Eq. (2),
it follows that
\[
P(1,i,t+1)=\left(1-\frac{1}{2t+N_{0}}\right)P(1,i,t).
\]
Since $P(1,t+1,t+1)=1$, one has
\begin{eqnarray*}
P(1,t+1)&=&\frac{1}{t+1+m_{0}}\sum\limits_{i=-m_{0},i\neq0}^{t+1}P(1,i,t+1)\\
&=&\frac{t+m_0}{t+1+m_0}\left(1-\frac{1}{2t+N_{0}}\right)P(1,t)+\frac{1}{t+1+m_0}.
\end{eqnarray*}
Then, by iteration,
\begin{eqnarray*}
P(1,t)&=&\prod\limits_{i=1}^{t-1}\left(1-\frac{1}{2i+N_{0}}\right)
\frac{i+m_0}{i+1+m_0}\left[P(1,1)+\sum\limits_{l=1}^{t-1}
\frac{\frac{1}{l+1+m_0}}{\prod\limits_{j=1}^{l}
\left(1-\frac{1}{2j+N_{0}}\right)\frac{j+m_0}{j+1+m_0}}\right]\\
&=&\frac{1}{t+m_0}\prod\limits_{i=1}^{t-1}\left(1-\frac{1}{2i+N_{0}}\right)
\left[(1+m_0)P(1,1)+\sum\limits_{l=1}^{t-1}\prod\limits_{j=1}^{l}
\left(1-\frac{1}{2j+N_{0}}\right)^{-1}\right].
\end{eqnarray*}

Next, let
\begin{eqnarray*}
x_n\triangleq&(1+m_0)P(1,1)
+\sum\limits_{l=1}^{n-1}\prod\limits_{j=1}^{l}
\left(1-\frac{1}{2j+N_{0}}\right)^{-1}
\end{eqnarray*}
and
\[
y_n\triangleq
(n+m_0)\prod\limits_{i=1}^{n-1}\left(1-\frac{1}{2i+N_{0}}\right)^{-1}>0.
\]
Thus, it follows that
\[
x_{n+1}-x_n=\prod\limits_{j=1}^{n}\left(1-\frac{1}{2j+N_{0}}\right)^{-1}
\]
and
\[
y_{n+1}-y_n=\frac{3n+N_{0}+m_0}{2n+N_{0}}\prod\limits_{i=1}^{n}
\left(1-\frac{1}{2i+N_{0}}\right)^{-1}>0.
\]
Since $y_n>0$ and $y_{n+1}-y_n>0$, $\{y_n\}$ is a strictly monotone
increasing nonnegative sequence, hence $y_n\rightarrow\infty$.
Moreover,
\[
\frac{x_{n+1}-x_n}{y_{n+1}-y_n}=\frac{2n+N_{0}}{3n+N_{0}+m_0}\
\rightarrow\ \frac{2}{3}~~ (n\rightarrow\infty).
\]

From Lemma 2, one has
\[
P(1)\triangleq\lim\limits_{t\rightarrow\infty}P(1,t)
=\lim\limits_{n\rightarrow\infty}\frac{x_n}{y_n}
=\lim\limits_{n\rightarrow\infty}\frac{x_{n+1}-x_n}{y_{n+1}-y_n}=\frac{2}{3}>0.
\]
This completes the proof.
\medskip

\textbf{Lemma 4} For $k>1$, if the limit
$\lim\limits_{t\rightarrow\infty}P(k-1,t)$ exists, then the limit
$\lim\limits_{t\rightarrow\infty}P(k,t)$ also exists and,
moreover,
\begin{eqnarray}
P(k)\triangleq\lim\limits_{t\rightarrow\infty}P(k,t)=\frac{k-1}{k+2}\,P(k-1)>0.
\end{eqnarray}

\textbf{Proof} First, observe that
\begin{eqnarray*}
P(k,t)&=&\frac{1}{t+m_0}\sum\limits_{i=-m_{0}i\neq0}^{t}P(k,i,t)\\
&=&\frac{1}{t+m_0}\sum\limits_{i=-m_{0}}^{-1}P(k,i,t)
+\frac{t}{t+m_0}\,\frac{1}{t}\,\sum\limits_{i=1}^{t}P(k,i,t).
\end{eqnarray*}

Next, denote
$\overline{P}(k,t)\triangleq\frac{1}{t}\sum\limits_{i=1}^{t}P(k,i,t)$.
One only needs to prove that the limit
$\lim\limits_{t\rightarrow\infty}\overline{P}(k,t)$ exists, which
will imply that the limit $\lim\limits_{t\rightarrow\infty}P(k,t)
=\lim\limits_{t\rightarrow\infty}\overline{P}(k,t)$ exists.

To show that the limit of $\overline{P}(k,t)$ exists as
$t\rightarrow\infty$, observe that $P(k,i,t)=0$ when $i>t+1-k$,
since in this case even if this vertex $i$ increases its degree by
1 each time, it cannot reach degree $k$. Then, it follows from
Lemma 1 that
\begin{eqnarray*}
\overline{P}(k,t)
&=&\frac{1}{t}\sum\limits_{i=1}^{t+1-k}P(k,i,t)\\
&=&\frac{1}{t}\sum\limits_{i=1}^{t+1-k}\sum\limits_{s=i+k-1}^{t}f(k,i,s)\,
   \prod\limits_{j=s}^{t-1}\left(1-\frac{k}{2j+N_{0}}\right)\\
&=&\frac{1}{t}\sum\limits_{i=1}^{t+1-k}\sum\limits_{s=i+k-1}^{t}
   P(k-1,i,s-1)\,\frac{k-1}{2(s-1)+N_{0}}\,\prod\limits_{j=s}^{t-1}
   \left(1-\frac{k}{2j+N_{0}}\right)\\
&=&\frac{1}{t}\sum\limits_{s=k}^{t}\sum\limits_{i=1}^{s+1-k}P(k-1,i,s-1)\,
   \frac{k-1}{2(s-1)+N_{0}}\,\prod\limits_{j=s}^{t-1}
   \left(1-\frac{k}{2j+N_{0}}\right)\\
&=&\frac{1}{t}\sum\limits_{s=k}^{t}\sum\limits_{i=1}^{s-1}
   P(k-1,i,s-1)\,\frac{k-1}{2(s-1)+N_{0}}\,\prod\limits_{j=s}^{t-1}
   \left(1-\frac{k}{2j+N_{0}}\right)\\
&=&\frac{1}{t}\sum\limits_{s=k}^{t}\overline{P}(k-1,s-1)\,
   \frac{(s-1)(k-1)}{2(s-1)+N_{0}}\,\prod\limits_{j=s}^{t-1}
   \left(1-\frac{k}{2j+N_{0}}\right)\\
&=&\frac{1}{t}\,\prod\limits_{i=k}^{t-1}\left(1-\frac{k}{2i+N_{0}}\right)
   \Bigg[\,\overline{P}(k-1,k-1)\frac{(k-1)^2}{2(k-1)+N_{0}}\\
& & +\sum\limits_{l=k}^{t-1}\overline{P}(k-1,l)\,
   \frac{l(k-1)}{2l+N_{0}}\,\prod\limits_{j=k}^{l}
   \left(1-\frac{k}{2j+N_{0}}\right)^{-1}\Bigg].
\end{eqnarray*}

Next, let
\begin{eqnarray*}
x_n&\triangleq&\overline{P}(k-1,k-1)\,\frac{(k-1)^2}{2(k-1)+N_{0}}\\
&&+\sum\limits_{l=k}^{n-1}\overline{P}(k-1,l)\,\frac{l(k-1)}{2l+N_{0}}
\prod\limits_{j=k}^{l}\left(1-\frac{k}{2j+N_{0}}\right)^{-1}
\end{eqnarray*}
and
\[
y_n\triangleq
n\prod\limits_{i=k}^{n-1}\left(1-\frac{k}{2i+N_{0}}\right)^{-1}>0
\rightarrow\infty.
\]
Obviously,
\[
x_{n+1}-x_n=\overline{P}(k-1,n)\,\frac{n(k-1)}{2n+N_{0}}\,
\prod\limits_{j=k}^{n}\left(1-\frac{k}{2j+N_{0}}\right)^{-1},
\]
and since
\begin{eqnarray*}
y_{n+1}-y_n
&=&\left[(n+1)-n\left(1-\frac{k}{2n+N_{0}}\right)\right]
   \prod\limits_{i=k}^{n}\left(1-\frac{k}{N_{0}+2i}\right)^{-1}\\
&=&\frac{(k+2)n+N_{0}}{2n+N_{0}}\,\prod\limits_{i=k}^{n}
   \left(1-\frac{k}{2i+N_{0}}\right)^{-1}>0,
\end{eqnarray*}
one has that $\{y_n\}$ is a strictly monotone increasing
nonnegative sequence, hence $y_n\rightarrow\infty$. Also, by
assumption,
\[
\frac{x_{n+1}-x_n}{y_{n+1}-y_n}=\frac{(k-1)n}{(k+2)n+N_{0}}\,
\overline{P}(k-1,n)\ \rightarrow\ \frac{k-1}{k+2}\,P(k-1)~~
(n\rightarrow\infty).
\]

Thus, it follows from Lemma 2 that
\[
\lim\limits_{t\rightarrow\infty}\overline{P}(k,t)
=\lim\limits_{n\rightarrow\infty}\frac{x_n}{y_n}
=\lim\limits_{n\rightarrow\infty}
\frac{x_{n+1}-x_n}{y_{n+1}-y_n}=\frac{k-1}{k+2}\,P(k-1)>0,
\]
therefore, $\lim\limits_{t\rightarrow\infty}P(k,t)$ exists, and Eq.
(6) is thus proved.
\medskip

\textbf{Theorem 1} The steady-state degree distribution of the BA
model with $m=1$ exists, and is given by
\begin{eqnarray}
P(k)=\frac{4}{k(k+1)(k+2)}\ \thicksim\ 4k^{-3}>0.
\end{eqnarray}

\textbf{Proof} By mathematical induction, it follows from Lemmas 3
and 4 that the steady-state degree distribution of the BA model
with $m=1$ exists. Then, solving Eq. (6) iteratively, one obtains
\[
P(k)=\frac{k-1}{k+2}\,P(k-1)
=\frac{k-1}{k+2}\,\frac{k-2}{k+1}\,\frac{k-3}{k}\,P(k-3).
\]
By continuing the process till $k=3+1$, one finally obtains
\[
P(k)=\frac{4}{k(k+1)(k+2)}\ \thicksim\ 4k^{-3}>0.
\]
\medskip

From Theorem 1, one can see that the degree distribution formula
of Krapivsky et al.$^{[5]}$ is exact, although the mathematical
proof there was not as rigorous as that given above.

\section{The BA model with $m\geq2$}

\ \ \ \ \textbf{Lemma 5} Under the $m\Pi$-hypothesis, the BA model
when $k>m$ satisfies
\begin{eqnarray}
f(k,i,s)&=&P(k-1,i,s-1)\,\frac{k-1}{2(s-1)+\frac{N_{0}}{m}},\\
P(k,i,t)&=&\sum\limits_{s=i+k-m}^t
           f(k,i,s)\,\prod\limits_{j=s}^{t-1}
           \left(1-\frac{k}{2j+\frac{N_{0}}{m}}\right).
\end{eqnarray}

\textbf{Proof} First, consider Eq. (8). The degree of vertex is
nondecreasing, and increasing at most by $1$ each time, according
to the construction of the BA model. Thus, it follows from the
Markovian properties that
\begin{eqnarray*}
f(k,i,s)&=&P\{k_i(s)=k,\ k_i(l)\neq k,\ l=1,2,\cdots,s-1\}\\
&=&P\{k_i(s)=k,\ k_i(s-1)=k-1,\ k_i(l)\neq k,\ l=1,2,\cdots,s-2\}\\
&=&P\{k_i(s)=k,\ k_i(s-1)=k-1\}\\
&=&P\{k_i(s-1)=k-1\}\,P\{k_i(s)=k|k_i(s-1)=k-1\}\\
&=&P(k-1,i,s-1)\,\frac{k-1}{2(s-1)+\frac{N_{0}}{m}}.
\end{eqnarray*}

Second, observe that the earliest time for the degree of vertex
$i$ to reach $k$ is at step $k+i-m$, and the latest time to do so
is at step $t$. After this vertex degree becomes $k$, it will not
increase any more. Thus, Eq. (9) is proved.
\medskip

\textbf{Lemma 6} Under the $m\Pi$-hypothesis, in the BA model the
limit $\lim\limits_{t\rightarrow\infty}P(m,t)$ exists and is
independent of the initial network; moreover,
\begin{eqnarray}
P(m)\triangleq\lim\limits_{t\rightarrow\infty}P(m,t)=\frac{2}{m+2}>0.
\end{eqnarray}

\textbf{Proof} From the construction of the BA model or (2), it
follows that
\[
P(m,i,t+1)=\left(1-\frac{m}{2t+\frac{N_{0}}{m}}\right)P(m,i,t).
\]
Since $P(m,t+1,t+1)=1$, one has
\begin{eqnarray*}
P(m,t+1)&=&\frac{1}{t+1+m_0}\,\sum\limits_{i=-m_{0},i\neq0}^{t+1}P(m,i,t+1)\\
&=&\frac{t+m_0}{t+1+m_0}\left(1-\frac{m}{2t+\frac{N_{0}}{m}}\right)P(m,t)
+\frac{1}{t+1+m_0}.
\end{eqnarray*}
Then, iterative calculation yields
\begin{eqnarray*}
P(m,t)&=&\prod\limits_{i=1}^{t-1}\left(1-\frac{m}{2i+\frac{N_{0}}{m}}\right)
\frac{i+m_0}{i+1+m_0}\left[\,P(m,1)+\sum\limits_{l=1}^{t-1}
\frac{\frac{1}{l+1+m_0}}{\prod\limits_{j=1}^{l}
\left(1-\frac{m}{2j+\frac{N_{0}}{m}}\right)\frac{j+m_0}{j+1+m_0}}\right]\\
&=&\frac{1}{t+m_0}\prod\limits_{i=1}^{t-1}
\left(1-\frac{m}{2i+\frac{N_{0}}{m}}\right)\left[(1+m_0)P(m,1)
+\sum\limits_{l=1}^{t-1}\prod\limits_{j=1}^{l}
\left(1-\frac{m}{2j+\frac{N_{0}}{m}}\right)^{-1}\right].
\end{eqnarray*}

Next, let
\[
x_n\triangleq(1+m_0)P(m,1)+\sum\limits_{l=1}^{n-1}
\prod\limits_{j=1}^{l}\left(1-\frac{m}{2j+\frac{N_{0}}{m}}\right)^{-1}
\]
and
\[
y_n\triangleq(n+m_0)\prod\limits_{i=1}^{n-1}
\left(1-\frac{m}{2i+\frac{N_{0}}{m}}\right)^{-1}>0.
\]
It follows that
\[
x_{n+1}-x_n
=\prod\limits_{j=1}^{n}\left(1-\frac{m}{2j+\frac{N_{0}}{m}}\right)^{-1}
\]
and
\[
y_{n+1}-y_n=\frac{(m+2)n+\frac{N_{0}}{m}+mm_0}{2n+\frac{N_{0}}{m}}\,
\prod\limits_{i=1}^{n}\left(1-\frac{m}{2i+\frac{N_{0}}{m}}\right)^{-1}>0.
\]
Consequently,
\[
\frac{x_{n+1}-x_n}{y_{n+1}-y_n}
=\frac{2n+\frac{N_{0}}{m}}{(m+2)n+\frac{N_{0}}{m}+mm_0}\
\rightarrow\ \frac{2}{m+2}~~ (n\rightarrow\infty).
\]
It follows from Lemma 2 that
\[
P(m)=\lim\limits_{t\rightarrow\infty}P(m,t)
=\lim\limits_{n\rightarrow\infty}\frac{x_n}{y_n}
=\lim\limits_{n\rightarrow\infty}\frac{x_{n+1}-x_n}{y_{n+1}-y_n}
=\frac{2}{m+2}>0.
\]
This completes the proof.
\medskip

\textbf{Lemma 7} Under the $m\Pi$-hypothesis, in the BA model with
$k>m$, if $\lim\limits_{t\rightarrow\infty}P(k-1,t)$ exists, then
$\lim\limits_{t\rightarrow\infty}P(k,t)$ exists and, moreover,
\begin{eqnarray}
P(k)\triangleq
\lim\limits_{t\rightarrow\infty}P(k,t)=\frac{k-1}{k+2}\,P(k-1)>0.
\end{eqnarray}

\textbf{Proof} First, observe that
\[
P(k,t)=\frac{1}{t+m_0}\sum\limits_{i=-m_{0},i\neq0}^{t}P(k,i,t)
=\frac{1}{t+m_0}\sum\limits_{i=-m_{0}}^{-1}P(k,i,t)
+\frac{t}{t+m_0}\frac{1}{t}\sum\limits_{i=1}^{t}P(k,i,t).
\]
Denote
$\overline{P}(k,t)\triangleq\frac{1}{t}\sum\limits_{i=1}^{t}P(k,i,t)$.
One only needs to prove that the limit
$\lim\limits_{t\rightarrow\infty}\overline{P}(k,t)$ exists, which
will imply that the limit $\lim\limits_{t\rightarrow\infty}P(k,t)
=\lim\limits_{t\rightarrow\infty}\overline{P}(k,t)$ exists.

To show that the limit of $\overline{P}(k,t)$ exists as
$t\rightarrow\infty$, observe that $P(k,i,t)=0$ when $i>t+m-k$.
Therefore, it follows from Lemma 5 that
\begin{eqnarray*}
\overline{P}(k,t)
&=&\frac{1}{t}\sum\limits_{i=1}^{t+m-k}P(k,i,t)\\
&=&\frac{1}{t}\sum\limits_{i=1}^{t+m-k}\sum\limits_{s=i+k-m}^{t}f(k,i,s)\,
   \prod\limits_{j=s}^{t-1}\left(1-\frac{k}{2j+\frac{N_{0}}{m}}\right)\\
&=&\frac{1}{t}\sum\limits_{i=1}^{t+m-k}\sum\limits_{s=i+k-m}^{t}
   P(k-1,i,s-1)\,\frac{k-1}{2(s-1)+\frac{N_{0}}{m}}\,\prod\limits_{j=s}^{t-1}
   \left(1-\frac{k}{2j+\frac{N_{0}}{m}}\right)\\
&=&\frac{1}{t}\sum\limits_{s=k-m+1}^{t}\sum\limits_{i=1}^{s+m-k}
   P(k-1,i,s-1)\,\frac{k-1}{2(s-1)+\frac{N_{0}}{m}}\,\prod\limits_{j=s}^{t-1}
   \left(1-\frac{k}{2j+\frac{N_{0}}{m}}\right)\\
&=&\frac{1}{t}\sum\limits_{s=k-m+1}^{t}\sum\limits_{i=1}^{s-1}
   P(k-1,i,s-1)\,\frac{k-1}{2(s-1)+\frac{N_{0}}{m}}\,\prod\limits_{j=s}^{t-1}
   \left(1-\frac{k}{2j+\frac{N_{0}}{m}}\right)\\
&=&\frac{1}{t}\sum\limits_{s=k-m+1}^{t}\overline{P}(k-1,s-1)\,
   \frac{(s-1)(k-1)}{2(s-1)+\frac{N_{0}}{m}}\,\prod\limits_{j=s}^{t-1}
   \left(1-\frac{k}{2j+\frac{N_{0}}{m}}\right)\\
&=&\frac{1}{t}\,\prod\limits_{i=k-m+1}^{t-1}
   \left(1-\frac{k}{2i+\frac{N_{0}}{m}}\right)
   \Bigg[\,\overline{P}(k-1,k-m)\frac{(k-1)(k-m)}{2(k-m)
   +\frac{N_{0}}{m}}\\
& &+\sum\limits_{l=k-m+1}^{t-1}\overline{P}(k-1,l)\,
   \frac{l(k-1)}{2l+\frac{N_{0}}{m}}\,\prod\limits_{j=k-m+1}^{l}
   \left(1-\frac{k}{2j+\frac{N_{0}}{m}}\right)^{-1}\Bigg]
\end{eqnarray*}

Next, let
\begin{eqnarray*}
x_n&\triangleq&\overline{P}(k-1,k-m)\,\frac{(k-1)(k-m)}{2(k-m)+\frac{N_{0}}{m}}\\
&&+\sum\limits_{l=k-m+1}^{n-1}\overline{P}(k-1,l)\,\frac{l(k-1)}{2l+\frac{N_{0}}{m}}
\prod\limits_{j=k-m+1}^{l}\left(1-\frac{k}{2j+\frac{N_{0}}{m}}\right)^{-1}
\end{eqnarray*}
and
\[
y_n\triangleq n\prod\limits_{i=k-m+1}^{n-1}
\left(1-\frac{k}{2i+\frac{N_{0}}{m}}\right)^{-1}>0.
\]
It follows that
\[
x_{n+1}-x_n=\overline{P}(k-1,n)\,\frac{n(k-1)}{2n+\frac{N_{0}}{m}}
\prod\limits_{j=k-m+1}^{n}\left(1-\frac{k}{2j+\frac{N_{0}}{m}}\right)^{-1}
\]
and
\begin{eqnarray*}
y_{n+1}-y_n
&=&\left[(n+1)-n\left(1-\frac{k}{2n+\frac{N_{0}}{m}}\right)\right]
   \prod\limits_{i=k-m+1}^{n}\left(1-\frac{k}{2i+\frac{N_{0}}{m}}\right)^{-1}\\
&=&\frac{(k+2)n+\frac{N_{0}}{m}}{2n+\frac{N_{0}}{m}}\prod\limits_{i=k-m+1}^{n}
   \left(1-\frac{k}{2i+\frac{N_{0}}{m}}\right)^{-1}>0.
\end{eqnarray*}
By assumption,
\[
\frac{x_{n+1}-x_n}{y_{n+1}-y_n}=\frac{(k-1)n}{(k+2)n+\frac{N_{0}}{m}}\,
\overline{P}(k-1,n)\rightarrow\frac{k-1}{k+2}P(k-1)~~
(n\rightarrow\infty).
\]
It then follows from Lemma 2 that
\[
\lim\limits_{t\rightarrow\infty}\overline{P}(k,t)
=\lim\limits_{n\rightarrow\infty}\frac{x_n}{y_n}
=\lim\limits_{n\rightarrow\infty}\frac{x_{n+1}-x_n}{y_{n+1}-y_n}
=\frac{k-1}{k+2}\,P(k-1)>0.
\]
Thus, $\lim\limits_{t\rightarrow\infty}P(k,t)$ exists and Eq. (11)
is proved.
\medskip

\textbf{Theorem 2} Under the $m\Pi$-hypothesis, the steady-state
degree distribution of the BA model with $m\geq2$ exists, and is
given by
\begin{eqnarray}
P(k)=\frac{2m(m+1)}{k(k+1)(k+2)}\ \thicksim\ 2m^2k^{-3}>0.
\end{eqnarray}

\textbf{Proof} By induction, it follows from Lemmas 6 and 7 that
the steady-state degree distribution of the BA model with $m\geq2$
exists. Equation (11) follows from iteration
\[
P(k)=\frac{k-1}{k+2}\,P(k-1)
=\frac{k-1}{k+2}\,\frac{k-2}{k+1}\,\frac{k-3}{k}\,P(k-3),
\]
till $k=3+m$. Thus, one obtains
\[
P(k)=\frac{2m(m+1)}{k(k+1)(k+2)}\ \thicksim\ 2m^2k^{-3}>0.
\]
\medskip

One can see that this degree distribution formula is consistent
with the formula obtained by Dorogovtsev et al.$^{[6]}$ and
Bollob\'as et al.$^{[9]}$, which allow multiple edges and loops.

\medskip
Finally, the authors thank Professor Yirong Liu for many helpful
discussions.

\vspace{.6cm}

Z. Hou$^1$:\ zthou@csu.edu.cn

X. Kong$^1$:\ kongxiangxing2008@163.com

D. Shi$^{1,2}$:\ shidh2001@263.net

G. Chen$^3$:\ gchen@ee.cityu.edu.hk